\documentclass[12pt,a4]{article}
\usepackage{amssymb}
\usepackage{latexsym}

\newtheorem{theorem}{Theorem}
\newtheorem{lemma}{Lemma}

\title{Groups of type $L_2(q)$ acting on polytopes}
\author{Dimitri Leemans\thanks{E-Mail: dleemans@ulb.ac.be}\\
Universit\'e Libre de Bruxelles\\
D\'epartement de Math\'ematiques - C.P.216\\
Boulevard du Triomphe\\
B-1050 Bruxelles\\[.06in]
{\small and}\\[.06in]
Egon Schulte\thanks{E-Mail: schulte@neu.edu}\\
Northeastern University\\
Department of Mathematics\\
360 Huntington Avenue\\
Boston, MA 02115, USA
}
\date{\today}
\begin{document}
\maketitle
\begin{abstract}
\noindent
We prove that if $G$ is a string C-group of rank 4 and $G\cong L_2(q)$ 
with $q$ a prime power, then $q$ must be 11 or 19. The polytopes arising are Gr\"unbaum's
11-cell of type $\{3,5,3\}$ for $L_2(11)$ and Coxeter's 57-cell of type $\{5,3,5\}$
for $L_2(19)$, each a locally projective regular $4$-polytope.
\end{abstract}

\section{Introduction}
\label{intro}

In this paper we determine the projective linear groups $L_2(q)$, $q$ a prime power, which 
occur as automorphism groups of abstract regular polytopes of rank $4$ (or higher). 

Given any type of group, the respective enumeration problem is generally difficult to solve due to the complexity of the subgroup lattice of the group at hand.  However, for small groups, an atlas of their regular polytopes is available (see Leemans \& Vauthier~\cite{LV2005}, Hartley~\cite{HA2005}). For example, the Suzuki group $Sz(8)$ admits regular polytopes, all of rank $3$, while
its automorphism group $Aut(Sz(8)) = Sz(8)\!:\!3$ does not. In fact, the following more
general result about almost simple groups of Suzuki type was obtained in Leemans~\cite{Lee2005a}. Let $G$ be a group such that $Sz(q) \leq G \leq Aut(Sz(q))$, where $q = 2^{2e+1}$ with $e$ a positive integer.  Then $G$ is a C-group if and only if $G = Sz(q)$. Moreover, if
$(G,\{\rho_0,\ldots,\rho_{n-1}\})$ is a string C-group, then $n=3$. Rephrased in terms of
regular polytopes this result says that, among such groups $G$, only the Suzuki group $Sz(q)$ itself can occur as automorphism group of a regular polytope of any rank, but does so only for regular polytopes of rank $3$, that is, for abstract regular polyhedra.

For groups of type $L_2(q)$ there is a wealth of interesting constructions of regular
polyhedra or maps on surfaces, but it is often assumed that $q$ is a prime (see, for
example,  McMullen~\cite{McM1991} or McMullen, Monson \& Weiss~\cite{MMW1993}). There is
also a considerable body of work available on representing $L_2(q)$ as a group of
orientation preserving mappings on surfaces, usually as the rotation (even) subgroup of the
automorphism group of a regular or chiral polyhedron or map (see, for example,  Glover \&
Sjerve~\cite{GS1985} or Conder~\cite{Co1990}, as well as Section~\ref{tworbits}). The
following result, established in Sjerve and Cherkassoff~\cite{SC94}, determines precisely
when a group $L_2(q)$ is the full automorphism group of a regular polyhedron, thereby
settling the existence problem when the rank is $3$. In fact, the group $L_2(q)$ may be
generated by three involutions, two of which commute, if and only if $q\neq 2, 3, 7$ or $9$. 

Regular polytopes of rank $5$ or higher cannot have an automorphism group of type $L_2(q)$
(see Leemans \& Vauthier~\cite{LV2005}). This, then, leaves the case of rank 4. The purpose of 
this article is to prove the following theorem, which was conjectured in \cite{LV2005} and
verified for $q\leq 61$. 

\begin{theorem}\label{ourthm}
If $L_2(q)$ is the automorphism group of a regular polytope of rank $4$, then $q=11$ or $19$.
\end{theorem}

In particular, it is found that there is just one regular polytope in each case. The {\em $11$-cell\/} with group $L_2(11)$ was constructed by Gr\"unbaum in \cite{Gr1977} by pasting eleven
hemi-icosahedra together (see also \cite{Cox1984}), and the {\em $57$-cell\/}  with
group $L_2(19)$ was discovered by Coxeter in \cite{Cox1982}. For more details see
Section~\ref{excep}.

\section{Preliminaries}
\label{def}

Abstract regular polytopes, string C-groups, and thin regular residually connected 
geometries with a linear diagram are basically the same mathematical objects. The link
between these objects has been discussed in, for instance, McMullen \& Schulte~\cite{MS2002}.
Here we take the viewpoint of string C-groups because it is the easiest and most efficient
one to describe abstract regular polytopes.

As defined in~\cite{MS2002}, a {\em C-group\/} is a group $G$ generated by pairwise 
distinct involutions $\rho_0,\ldots,\rho_{n-1}$ which satisfy the following property,
called the {\it intersection property\/}:
\[\forall J, K \subseteq \{0,\ldots,n-1\}, \;\; 
\langle \rho_j \mid j \in J\rangle \cap \langle \rho_j \mid j \in K\rangle = 
\langle\rho_j \mid j \in J\cap K\rangle . \]
We call $n$ the {\em rank\/} of $G$.

A C-group $(G, \{\rho_0,\ldots, \rho_{n-1}\})$, or simply $G$, is a {\em string\/} 
C-group if its generators satisfy the following relations:
\[ (\rho_j\rho_k)^2 = 1 \;\; 
\forall j, k \in \{0,\ldots n-1\} \mbox{ with} \mid j-k\mid \geq 2.\]

Each string C-group $G$ then determines (uniquely) a regular $n$-polytope $\cal P$ with 
automorphism group $G$. The {\em $i$-faces} of $\cal P$ are the right cosets of the
distinguished subgroup $G_i := \langle \rho_{j} \mid j \neq i \rangle$
for each $i = 0,1,\ldots, n-1$, and two faces are incident just when they intersect as
cosets.  Formally,  we also adjoin two copies of $G$ itself, as the (unique) $(-1)$- and
$n$-faces of $\cal P$.  Conversely, the automorphism group of a regular $n$-polytope is a
C-group, whose generators $\rho_j$ map a fixed, or {\em base\/}, flag $\Phi$ of $\cal P$ to
the $j$-adjacent flag $\Phi^{j}$ (differing from $\Phi$ in the $j$-face).

Recall from \cite{MS2002} that $\{{\cal P}_1,{\cal P}_2\}$ denotes the universal regular
polytope (if it exists) with facets ${\cal P}_1$ and vertex-figures ${\cal P}_2$. This
covers every regular polytope with facets ${\cal P}_1$ and vertex-figures ${\cal P}_2$.

\section{The subgroup structure of $L_2(q)$}
\label{pslsub}

We frequently require properties of the subgroup lattice of $L_2(q)$. The problem of
determining all the subgroups of $L_2(q)$ has been completely solved by Dickson
\cite{Dic58}. Another proof can be found in \cite{Hup67}.

\begin{theorem}\label{sub}

The group $L_2(q)$ of order $\frac{q  (q^2 - 1)}{(2,q-1)}$, where $q=p^r$ with $p$ a
prime, contains only:

\begin{enumerate}
\item $q+ 1$ conjugate elementary abelian subgroups of order $q$, denoted by $E_q$.

\item $\frac{q(q\mp 1)}{2}$ conjugate cyclic subgroups of order $d$,  denoted by $d$, for all divisors $d$ of $\frac{(q\pm 1)}{(2,q-1)}$.

\item $\frac{q(q^2 -1)}{2d (2,q-1)}$ dihedral groups of order $2d$, denoted by $D_{2d}$,
for all divisors $d$ of $\frac{(q\pm 1)}{(2,q-1)}$ with $d>2$. The number of conjugacy
classes of these subgroups is one if $\frac{(q\pm 1)}{d (2,q-1)}$ is odd, and two if
it is even.

\item For $q$ odd, $\frac{q(q^2 -1)}{12 (2,q-1)}$ dihedral groups of order 4, denoted
by $2^2$. The number of conjugacy classes of these groups is one if  $q \equiv \pm 3 (
8)$ and two if $q \equiv \pm 1 ( 8)$. For $q$ even, the groups $2^2$ are listed under
family~5.

\item $\frac{(2,1,1) (p^k -1)(p^r -1)(p^r -p) \cdot \cdot \cdot (p^r -p^{s-1})}
{(q-1) (p^s -1)(p^s -p) \cdot \cdot \cdot (p^s -p^{s-1})}$ sets, each of
$\frac {q^2 -1} {(2,1,1) (p^k -1)}$ conjugate elementary abelian subgroups of order 
$p^s$, denoted by $E_{p^s}$, for all natural number $s$ such that $1\leq s \leq r-1$,
where $k=(r,s)$ and $(2,1,1)$ is defined as 2, 1 or 1 according as $p>2$ and 
$\frac{r}{k}$ is even, $p>2$ and $ \frac {r}{k}$ is odd, or $p=2$.

\item $\frac{(2,1,1) (p^k -1)(p^r -p) \cdot \cdot \cdot (p^r -p^{s-1})}
{p^{r-s}(q-1) (p^s -1)(p^s -p) \cdot \cdot \cdot (p^s -p^{s-1})}$
sets of $\frac {(q^2 -1) p^{r-s}} {(2,1,1) (p^k -1)} $ subgroups $E_{p^s}\!:\!h$,
each a semidirect product of an elementary abelian group $E_{p^s}$ and a cyclic 
group of order $h$, for all natural numbers $s$ such that $1\leq s \leq r$ and
all divisors $h$ of $\frac {p^k -1}{(2,1,1)}$, where again $k=(r,s)$ and $(2,1,1)$ is defined as 2,
1 or 1 according as $p>2$ and  $\frac {r}{k}$ is even, $p>2$ and $\frac {r}{k}$ is
odd, or $p=2$.

\item For $q$ odd or $q=4^m$,\  $\frac{q(q^2 -1)}{12 (2,q-1)}$ alternating groups $A_4$,
of order 12. The number of conjugacy classes of these groups is one if 
$q \equiv \pm 3 (8)$ and two if $q \equiv \pm 1 ( 8)$.

\item For $q \equiv \pm 1 ( 8)$, two conjugacy classes of $\frac{q(q^2 -1)}{24 (2,q-1)}$
symmetric groups $S_4$, of order 24.

\item For $q \equiv \pm 1 ( 5)$, two conjugacy classes of $\frac{q(q^2 -1)}{60 (2,q-1)}$
alternating groups $A_5$, of order 60; and for $q=4^m$, one conjugacy class of
$\frac{q (q^2 -1)}{60 (2,q-1)}$ alternating groups $A_5$. For $q \equiv 0 ( 5)$, the
groups $A_5$ are listed under family~10.

\item $\frac{q(q^2 -1)}{ p^w (p^{2w}-1)}$ groups $L_2(p^w)$, for all divisors $w$ of
$r$. The number of conjugacy classes of these groups is two, one or one according as $p>2$
and $\frac {r}{w}$ is even, $p>2$ and $ \frac {r}{w}$ is odd, or $p=2$.

\item Two conjugacy classes of $\frac{q(q^2 -1)}{2 p^w (p^{2w}-1)}$ groups $PGL_2(p^w)$,
for all $w$ such that $2w$ is a divisor of $r$.
\end{enumerate}
\end{theorem}

Observe that when $q$ is even, family~11 of Theorem~\ref{sub} is a subfamily of family~10.

\section{$L_2(q)$ acting flag-transitively}
\label{pslacts}

In this section, we assume that $G$ is a group isomorphic to $L_2(q)$, with $q=p^r$, 
$p$ a prime and $r$ a positive integer. Moreover, we assume that $(G, \{\rho_0,\ldots,
\rho_{3}\})$ is a string C-group of type $\{t,l,s\}$ (i.e. the orders of
$\rho_0\rho_1$, $\rho_1\rho_2$ and $\rho_2\rho_3$ are $t$, $l$ and $s$, respectively). 
Clearly, $t,l,s\geq 3$, since $G$ is not a direct product of two non-trivial groups. As
before we set 
\[ G_i = \langle \rho_j \mid j \in \{0,\ldots, 3\}\backslash\{i\}\rangle, \] 
for $i=0,\ldots,3$. Our aim is to prove that $q$ must be $11$ or $19$.

We say that a subgroup $H$ of $G$ is an (irreducible) {\em rank $3$ subgroup\/} of $G$ if $H$ is a rank 3 string C-group with a connected Coxeter diagram.

By Theorem~\ref{sub}, the rank 3 subgroups of $G$ are isomorphic to $S_4$, $A_5$, or
$L_2(q')$ or $PGL_2(q')$ for some $q'$. These are the only possible types of subgroups for
$G_0$ and $G_3$.

We begin with a sequence of lemmas aimed at eliminating $L_2(q')$ and $PGL_2(q')$ as
possibilities.

\begin{lemma}\label{lemma1}
The prime $p$ must be odd.
\end{lemma}
{\bf Proof:}
The subgroup $G_2 = \langle \rho_0,\rho_1,\rho_3\rangle$ is a subgroup of the form 
$2\times D_{2t}$, with $t$ as above. By Theorem~\ref{sub}, this subgroup must be contained in
a maximal subgroup of dihedral type of $G$. Hence $q$ must be odd.
\hfill $\Box$
\bigskip

From now on we may assume that $p$ is {\em odd}.

\begin{lemma}\label{lemma2}
The orders $t$ of $\rho_{0}\rho_{1}$ and $s$ of $\rho_{2}\rho_{3}$ must be odd.
\end{lemma}
{\bf Proof:}
Let us prove, without loss of generality, that $t$ must be odd. Suppose $t$ is even.
Then $G_2\cong 2^2\times D_t$. Inspecting the list of subgroups of $L_2(q)$ given by
Theorem~\ref{sub}, we readily see that this cannot occur. Therefore, $t$ (and $s$) must
be odd.
\hfill $\Box$

\begin{lemma}\label{lemma3}
Let $H$ and $K$ be two subgroups of type $L_2(q')$ in $L_2(q)$, with 
$q'^{\,m} = q$ for some positive integer $m$. Then $H\cap K$ cannot be a dihedral group
$D_{2k}$ with $k>2$ (and $k$ a divisor of $\frac{(q'\pm 1)}{2}$).
\end{lemma}
{\bf Proof:}
Let $k>2$, and let $k$ be a divisor of $\frac{q'\pm 1}{2}$. By Theorem~\ref{sub}, we
know that 
\begin{itemize}
\item in $L_2(q)$, there are $\frac{q(q^2-1)}{q'(q'^2-1)}$ subgroups isomorphic to 
$L_2(q')$;
\item in $L_2(q)$, there are $\frac{q(q^2-1)}{4k}$ subgroups isomorphic to $D_{2k}$;
\item in $L_2(q')$, there are $\frac{q'(q'^2-1)}{4k}$ subgroups isomorphic to
$D_{2k}$.
\end{itemize}
Let $n := \frac{(q'-1)}{2}$ if $k\mid \frac{(q'-1)}{2}$ and $n:=\frac{(q'+1)}{2}$ if $k\mid \frac{(q'+1)}{2}$.
By Theorem~\ref{sub}, there are $\frac{q(q^2-1)}{4n}$ subgroups $D_{2n}$ in $L_2(q)$. Each subgroup $D_{2n}$ contains $\frac{n}{k}$ subgroups $D_{2k}$.
Therefore, each subgroup $D_{2k}$ is contained in exactly one subgroup $D_{2n}$.
The same kind of arguments show that each $D_{2n}$ is contained in exactly one $L_2(q')$.
Therefore, each subgroup $D_{2k}$ of $L_2(q)$ (with $k>2$ and $k$ a divisor
of $\frac{(q'\pm 1)}{2}$) is contained in a subgroup $L_2(q')$ of $L_2(q)$, and the number of subgroups $L_2(q')$ containing a given subgroup $D_{2k}$ is precisely one. Now the lemma follows.
\hfill $\Box$

\begin{lemma}\label{lemma4}
The subgroups $G_0$ and $G_3$ of $G$ cannot be isomorphic to $L_2(q')$, with
$q'^{\,m} = q$ for some positive integer $m$.
\end{lemma}
{\bf Proof:}
Suppose, without loss of generality, that $G_3 \cong L_2(q')$. Then $\rho_3$
conjugates two subgroups $L_2(q')$ of $G$ whose intersection contains a dihedral group of order
$2t$; in terms of the underlying polytope, of which $G$ is the automorphism group, these
two subgroups are the stabilizers of the two facets which share the $2$-face in the base flag.
Then this intersection itself, being a subgroup of $L_2(q)$, must be a dihedral group.
However, this  contradicts Lemma~\ref{lemma3}.~\hfill $\Box$

\begin{lemma}\label{lemma5}
Let $H$ and $K$ be two subgroups of type $L_2(q')$ in $L_2(q)$, with
$q'^{\,m} = q$ for some positive integer $m$. If $H\cap K$ contains a cyclic group $Z_k$,
$k>2$, whose normalizer $N_{L_2(q)}(Z_k)$ in $L_2(q)$ is a maximal subgroup of
$L_2(q)$ of dihedral type, then $H\cap K\cong Z_\frac{(q'-1)}{2}$ when
$k\mid \frac{(q'-1)}{2}$ and $H\cap K \cong Z_\frac{(q'+1)}{2}$ when $k\mid
\frac{(q'+1)}{2}$.
\end{lemma}
{\bf Proof:}
As before, let $n := \frac{(q'-1)}{2}$ if $k\mid \frac{(q'-1)}{2}$ and $n:=\frac{(q'+1)}{2}$ if 
$k\mid \frac{(q'+1)}{2}$. For the normalizers in $H$ and $K$ we have $N_H(Z_k) = N_H(Z_n) = D_{2n}$ and $N_K(Z_k) = N_K(Z'_{n}) = D'_{2n}$, where $D_{2n}$ and $D'_{2n}$ are dihedral subgroups of $H$ and $K$ of order $2n$, respectively, contained in the (dihedral) normalizer 
$D:=N_{L_2(q)}(Z_k)$ (of order $q\pm 1$), and $Z_{n}$ and $Z'_{n}$ are their cyclic subgroups of order $n$. However, $Z_{n}$ and $Z'_{n}$ are cyclic subgroups of $D$ of the same order, hence $Z_{n}=Z'_{n}\leq H\cap K$. By Lemma~\ref{lemma3}, $H\cap K$ cannot be dihedral, so 
$H\cap K = Z_{n}$. 
\hfill $\Box$

\begin{lemma}\label{lemma7}
The subgroups $G_0$ and $G_3$ of $G$ cannot be isomorphic to $PGL_2(q')$, with
$q'^{\,2m}=q$ for some positive integer $m$.
\end{lemma}
{\bf Proof:}
Suppose, without loss of generality, that $G_3 \cong PGL_2(q')$. Then, by the same argument
as in Lemma~\ref{lemma4}, $\rho_3$ conjugates two subgroups $PGL_2(q')$, namely $G_3$ and $G_{3}^{\rho_3}$, and their intersection $D$ must contain the dihedral group 
$D_{2t} := \langle\rho_{0},\rho_{1}\rangle$ of order $2t$. 

First observe that if $m > 1$, there exists a subgroup $L\cong L_2(q'^{\,2})$ such that $G_3
< L < G$. The element $\rho_3$ fixes $D$ by conjugation, namely $D^{\rho_3} = D$, and hence 
$D\leq L\cap L^{\rho_3}$. Therefore, $L\cap L^{\rho_3}$ itself must be a dihedral group,
which contradicts Lemma~\ref{lemma3}. Hence we may assume that $m=1$.

Each of $G_3$ and $G_3^{\rho_{3}}$ contains a subgroup of index 2 isomorphic to $L_2(q')$. Since the intersection of these two subgroups $H$ and $K\,(=H^{\rho_3})$ cannot be dihedral by Lemma~3, it must coincide with a cyclic subgroup $Z_k$ of $D$ of some order $k$.  Then the normalizer of $Z_k$ in $G$ is a maximal subgroup of dihedral type, since it contains $D$ (and $\rho_3$).  By Lemma~\ref{lemma5}, we then know that $k$ must be $\frac{(q'- 1)}{2}$ or $\frac{(q'+ 1)}{2}$; that is, $k=n$, in our previous notation. On the other hand, since $t$ is odd by Lemma~\ref{lemma2} and the square of every element in $D$ (and hence of $\rho_0\rho_1$) is necessarily in $H\cap K$ (the index of $L_2(q')$ in $PGL_2(q')$ is $2$), we also have $t\mid n$, that is, $t \mid n=k$. 

We claim that $D$ is a dihedral group of order $4n$. First recall that the maximal dihedral subgroups of $PGL_2(q')$ are of order $2(q'\pm 1)$ (see Moore~\cite{Moo1904}), so certainly $D$ is of order $2n$ or $4n$. Next observe that every subgroup $D_{2n}$ of $G$ is contained in a unique subgroup $D_{4n}$ of $G$. This can be seen as follows.  Clearly, every dihedral subgroup $D_{4n}$ of $G$ contains exactly two dihedral subgroups $D_{2n}$.  Moreover, by item~(3) of Theorem~\ref{sub}, the number of subgroups $D_{2n}$ in $G$ is exactly twice the number of subgroups $D_{4n}$ in $G$. Hence, since every subgroup $D_{2n}$ of $G$ actually is contained in a subgroup $D_{4n}$ of $G$, the latter must necessarily be uniquely determined, proving our claim. But now we can argue as follows. Assume that $D$ is only of order $2n$. Then $D$, viewed as a dihedral subgroup of $G_3$ of order $2n$, is contained in a maximal dihedral subgroup $D'$ of $G_3$ of order $4n$. Similarly, $D$, viewed as a dihedral subgroup of $G_{3}^{\rho_3}$, is also contained in a maximal dihedral subgroup $D''$ of $G_{3}^{\rho_3}$ of order $4n$.  Hence, by uniqueness in $G$, 
\[ D\leq D'=D'' \leq G_{3} \cap G_{3}^{\rho_3} = D ; \]
that is, $D=D'$, of order $4n$.  This is a contradiction, so $D$ must be of order $4n$.

Next we proceed by constructing a dihedral subgroup $E$ in $D$ of order $2t$ which is contained in $H$ and $K$. This, then, forces $H\cap K$ to be dihedral and once again provides a  contradiction to Lemma~\ref{lemma3}, thereby completing the proof.

Let $\tau_0,\tau_1$ be a pair of involutory generators of $D$ with $(\tau_0\tau_1)^{2n}=1$. Since 
$D$ contains $\langle\rho_{0},\rho_{1}\rangle = D_{2t}$, we may assume that $\tau_{0}=\rho_{0}$.
Then $(\tau_0\tau_1)^{2}\in H\cap K$ but $\tau_0\tau_1\not\in H\cap K$. Hence $\tau_0\tau_1\not\in H$ or $\tau_0\tau_1\not\in K$. On the other hand, $\rho_{0}\rho_{1} \in H \cap K$. However, $\rho_{0}\not\in H$ and $\rho_{0}\not\in K$; in fact, if $\rho_{0}$ is in $H$ or $K$, respectively, then 
$\rho_{0}=\rho_{0}^{\rho_{3}}$ is in $K\,(=H^{\rho_3})$ or $H\,(=K^{\rho_3})$ and hence $\tau_{0}=\rho_{0}\in H\cap K = \langle (\tau_0\tau_1)^{2} \rangle$, which is impossible. Now suppose, without loss of generality, that $\tau_0\tau_1\not\in H$. Then, since $\tau_{0}=\rho_{0}\not\in H$ and $H$ has index $2$ in $G_3$, we must have $\tau_{1}\in H$. It follows that the subgroup $E := \langle \rho_{0}\rho_{1},\tau_{1} \rangle$ of $D$ must be a dihedral group of order $2t$ contained in $H$. We claim that $E$ is also contained in $K$ and hence in $H\cap K$.  For the proof we need to verify that 
$\tau_{1}\in K$.  Assume to the contrary that $\tau_{1}\not\in K$. Then, since $\tau_{0}=\rho_{0}\not\in K$ and $K$ has index $2$ in $G_{3}^{\rho_3}$, we must have $\tau_{0}\tau_{1}\in K$. However, then $K$ would contain the cyclic subgroup $\langle\tau_{0}\tau_{1}\rangle$ of order $2n$, contradicting the fact that the maximal cyclic subgroups of a group $L_2(q')$ are of order $\frac{(q'\pm 1)}{2}$. It follows that $\tau_{1}\in K$ and hence $E\leq K$, as required. This completes the proof.
\hfill $\Box$
\medskip

We finally have all the tools to prove the following theorem, which in turn implies Theorem~\ref{ourthm}.

\begin{theorem}
Let $(G, \{\rho_0,\ldots, \rho_{3}\})$ be a string C-group.
Suppose that $G\cong L_2(q)$. Then $q=11$ or $19$.
\end{theorem}
{\bf Proof:}
Lemmas~\ref{lemma3} and~\ref{lemma7} reduce the possible subgroups for $G_0$ and $G_3$ to
only two kinds, $S_4$ and $A_5$.

The only rank 3 polytopes with group $S_4$ are $\{3,3\} \,(=\{3,3\}_4)$, $\{3,4\}_3$ and
$\{4,3\}_3$, and those with group $A_5$ are $\{3,5\}_5$, $\{5,3\}_5$ and $\{5,5\}_3$.
Recall here from \cite{CM1980} that $\{m,n\}_k$ is obtained from the regular tessellation
$\{m,n\}$ by identifying any two vertices that are separated by $k$ steps along a Petrie
polygon of $\{m,n\}$. Its group $\langle\tau_0,\tau_1,\tau_2\rangle$ has a presentation consisting of the standard Coxeter type relations for $\{m,n\}$ and the single extra relation 
$(\tau_{0}\tau_{1}\tau_{2})^{k} = 1$.)  

We can now check which pairs of polyhedra can be combined to form the facets and vertex-figures, 
respectively, of a regular rank 4 polytope. Table~\ref{table1} gives the possible
combinations and the structure of the corresponding ``universal" groups; these groups are
obtained by taking as defining relations just those of the facet group and vertex-figure
group as well as $(\rho_{0}\rho_{3})^2 = 1$.  Only one from a pair of dual combinations is
listed, since dual combinations yield the same groups (with the orders of the generators reversed). The results in this table can easily be obtained using a Computational Algebra package like {\sc
Magma}~\cite{BCP97} (or, if necessary, by hand). Finally, by inspection we readily see that
the only possibilities for $(G,\{\rho_0,\ldots,\rho_{3}\})$ to be a string C-group of rank
$4$ occur when $q=11$ or $q=19$. 
\hfill $\Box$

\begin{table}
\begin{center}
\begin{tabular}{||c|c|c|c||}
\hline
Facet&Vertex-figure&Order of $G$&Structure of $G$\\
\hline
$\{5,5\}_3$&$\{5,5\}_3$&1&\\
$\{5,5\}_3$&$\{5,3\}_5$&1&\\
$\{5,3\}_5$&$\{3,5\}_5$&3420&$L_2(19)$\\
$\{5,3\}_5$&$\{3,4\}_3$&60&$A_5 \cong L_2(5)$\\
$\{5,3\}_5$&$\{3,3\}_4$&1&\\
$\{4,3\}_3$&$\{3,4\}_3$&96&$2^4\! : \!S_3$ \\
$\{4,3\}_3$&$\{3,3\}_4$&24&$S_4$ \\
$\{3,5\}_5$&$\{5,3\}_5$&660&$L_2(11)$\\
$\{3,4\}_3$&$\{4,3\}_3$&1&\\
$\{3,3\}_4$&$\{3,3\}_4$&120&$S_5$\\
\hline
\end{tabular}
\caption{Combinations of rank 3 polytopes}\label{table1}
\end{center}
\end{table}

Note that the groups occurring in rows 1, 2, 5 and 9 of Table~\ref{table1} are trivial. In
row~4, the group actually is a group $L_2(q)$ but is too small to be a C-group of rank 4.
In row~7 we also do not have a C-group of rank 4. Finally, in row~6 we obtain a C-group of
rank 4 isomorphic to $2^4\! : \!S_3$, namely the group of the universal locally projective
regular $4$-polytope $\{\{4,3\}_{3},\{3,4\}_{3}\}$.   

\section{The 11-cell and 57-cell}
\label{excep}

In terms of regular polytopes our main theorem can be rephrased as follows. 

\begin{theorem}
The only regular polytopes of rank $4$ with automorphism groups of type $L_2(q)$ are the
$11$-cell $\{\{3,5\}_{5},\{5,3\}_{5}\}$ with group $L_2(11)$ of order $660$, and the
$57$-cell $\{\{5,3\}_{5},\{3,5\}_{5}\}$ with group $L_2(19)$ of order $3420$.
\end{theorem}

The two polytopes of the theorem are self-dual and locally projective (see \cite{MS2002}); their 
facets and vertex-figures are regular maps in the projective plane. The $11$-cell has $11$
hemi-icosahedral facets, $11$ vertices with hemi-dodecahedral vertex-figures, $55$ edges,
and $55$ triangular $2$-faces. (The hemi-icosahedron and hemi-dodecahedron, respectively, are
obtained from the icosahedron or dodecahedron by identifying pairs of antipodal vertices,
edges and $2$-faces.) The edge-graph of the $11$-cell is the complete graph on $11$
vertices. The $57$-cell has $57$ hemi-dodecahedral facets, $57$ vertices with
hemi-icosahedral vertex-figures, $171$ edges, and $171$ pentagonal $2$-faces. Each polytope
is universal among the regular polytopes with the same kind of facets and vertex-figures. 
Moreover, in terms of their basic generators $\rho_0,\ldots,\rho_3$, each of the groups
$L_2(11)$ and $L_2(19)$ has a presentation consisting of the standard Coxeter type
relations for the $3$-dimensional regular hyperbolic honeycomb $\{3,5,3\}$ or $\{5,3,5\}$,
respectively, as well as the two extra relations
\[ (\rho_{0}\rho_{1}\rho_{2})^{5} = (\rho_{1}\rho_{2}\rho_{3})^{5} = 1 \]
(the same two extra relations in each case).

It is interesting to note the effect of dropping one, the first (say), of the two extra relations just mentioned. For the type $\{5,3,5\}$ we then obtain the universal regular $4$-polytope 
$\{\{5,3\},\{3,5\}_{5}\}$ with automorphism group 
\[ J_{1} \times L_2(19) , \] 
where $J_1$ denotes the first Janko group, a sporadic simple group of order 175560 (see 
\cite{HL2004}). This polytope covers the $57$-cell.  

In \cite{Gr1977}, Gr\"unbaum showed that there is no polytope of type $\{3,5,3\}$ with icosahedral facets and hemi-dodecahedral vertex figures (see also \cite{HA2004}). Hence, in that case, dropping one of the two extra relations still gives $L_2(11)$.

Recall that the even (or rotation) subgroup of a C-group consists of the elements which can be 
expressed as products of an even number of generators $\rho_i$; it has index at most $2$ in
the full group.  Note that, when considered as C-groups, $L_2(11)$ and $L_2(19)$
coincide with their even subgroups, because they are simple groups. 

\section{$L_2(q)$ acting with two flag orbits}
\label{tworbits}

Our main theorem can be rephrased as saying that there are just two abstract polytopes of rank 
$4$ on which a group of type $L_2(q)$ admits a flag-transitive action as a group of
automorphisms. By contrast, there is a wealth of abstract polytopes of rank $4$ on which a
group $L_2(q)$ acts with precisely two flag orbits. In fact, this is already true in rank
$3$, the most prominent example being $L_2(7)$, acting as the even subgroup of the
automorphism group $PGL_2(7)$ of Klein's regular map $\{3,7\}_8$ of genus $3$ (see
\cite{CM1980}, as well as \cite{GS1985})). 

The known examples of polytopes $\cal P$ with two flag orbits generally have the property
that adjacent flags are in distinct orbits. This can arise in one of two ways. Either the
polytope $\cal P$ is regular and $L_2(q)$ is the even subgroup of its automorphism group
(as in the case of the Klein map), or $\cal P$ is chiral and $L_2(q)$ is its full
automorphism group. (Recall here that an abstract polytope is {\em chiral\/} if its automorphism
group has two orbits on the flags, such that adjacent flags are in distinct orbits.)

One type of construction of examples of rank $4$ begins with a $3$-dimensional regular hyperbolic honeycomb and a faithful representation of its symmetry group as a group of complex M\"obius transformations, generated by the inversions in four circles cutting one another at the same angles as the corresponding reflection planes in hyperbolic space (see \cite{SW1994}). For example, for the regular honeycomb $\{4,4,3\}$, the even subgroup of its symmetry group is isomorphic to $L_2(\mathbb{Z}[i]) \!:\! 2$, where $\mathbb{Z}[i]$ is the ring of Gaussian integers. Then interesting finite regular or chiral polytopes $\cal P$ of rank $4$ can be obtained by
modular reduction of this group. In each case, the resulting group is of type $L_2(R)$,
with $R$ a finite ring, and acts as a group of automorphisms of $\cal P$ with two flag
orbits. In certain cases, $R$ is a field.  

For example, if $p$ is a prime with $p \equiv 3 (4)$, then $\cal P$ is regular, the even 
subgroup is $L_2(p^2)$,  and $\cal P$ has toroidal facets $\{4,4\}_{(p,0)}$ and cubical
vertex-figures $\{4,3\}$ (see \cite[p.238]{SW1994}). Similarly, if $p$ is a prime with $p
\equiv 1 (8)$, then $\cal P$ is chiral, the automorphism group is $L_2(p)$, and $\cal P$
has toroidal facets $\{4,4\}_{(b,c)}$, with $b^2 + c^2 = p$, and again cubical vertex-figures
$\{4,3\}$ (see \cite[p.239]{SW1994}). There are similar such results for other Schl\"afli
symbols (see also \cite{NS1995}).

Yet more examples of polytopes with two flag orbits can be found in, for example, \cite{MWei1990}.

{\bf Acknowledgement}
This research was accomplished while the first author was visiting the second author at
Northeastern University. The first author gratefully acknowledges financial support from the
Belgian National Fund for Scientific Research, and the ``Fondation Agathon De Potter" from
the ``Classe des Sciences" of the ``Acad\'emie Royale de Belgique". The second author was
supported by NSA-grant H98230-05-1-0027.

\bibliographystyle{plain}

\end{document}